\documentclass[12pt]{amsart}
\usepackage{amssymb}

\newtheorem{fact}{Fact}[section]

\newtheorem{theorem}[fact]{Theorem}
\newtheorem{definition}[fact]{Definition}
\newtheorem{example}[fact]{Example}
\newtheorem{rremark}[fact]{Remark}
\newtheorem{proposition}[fact]{Proposition}
\newtheorem{corollary}[fact]{Corollary}
\newenvironment{remark}{\begin{rremark} \rm}{\end{rremark}}

\def\mapright#1{\smash{\mathop{\longrightarrow}\limits^{#1}}}
\def\mapleft#1{\smash{\mathop{\longleftarrow}\limits^{#1}}}
\def\mapdown#1{\Big\downarrow\rlap{$\vcenter{\hbox{$\scriptstyle#1$}}$}}

\DeclareMathOperator{\ts}{Ts} \DeclareMathOperator{\tp}{Tp}
\def\ss{{}^{\sharp}}
\def\fl#1{{}^{\hskip.08cm\flat[#1]}}
\DeclareMathOperator{\E}{\mathcal E}
\DeclareMathOperator{\C}{\mathbf C} \DeclareMathOperator{\Z}{\mathbf Z}
\DeclareMathOperator{\K}{\mathcal K}
\DeclareMathOperator{\A}{\mathcal A}
\DeclareMathOperator{\Q}{\mathbf Q}
\DeclareMathOperator{\Hom}{Hom}

\title{On the structure of Thom polynomials of singularities}

\author{L. M. Feh\'er}
\address{Department of Analysis, Eotvos University Budapest, Hungary}
\email{lfeher@renyi.hu}

\author{R. Rim\'anyi}
\address{Department of Mathematics, University of North Carolina at Chapel Hill, USA}
\email{rimanyi@email.unc.edu}

\begin{document}
\thanks{\noindent Supported by OTKA T046365MAT, and
        NSF grant DMS-0405723 (2nd author) \\
Keywords: Thom polynomial, equivariant maps\\
AMS Subject classification 32S20}

\begin{abstract} Thom polynomials of singularities express the cohomology classes dual to singularity submanifolds. A stabilization
property of Thom polynomials is known classically, namely that trivial unfolding does not change the Thom polynomial. In this paper we show that this is a special case of a product rule. The product rule enables us to calculate the Thom polynomials of singularities if we know the Thom polynomial of the product singularity. As a special case of the product rule we  define a formal power series (Thom series, $\ts_Q$) associated with a commutative, complex, finite dimensional local algebra $Q$, such that the Thom polynomial of {\em every} singularity with local algebra $Q$ can be recovered from $\ts_Q$.
\end{abstract}

\maketitle

\section{Introduction}

For a holomorphic map $f:N\to M$ between complex manifolds, and a
singularity $\eta$, one can consider the points in $N$ where the map
has singularity $\eta$. In \cite{thom_orig} Ren\'e Thom noticed that
the cohomology class represented by this set can be calculated by a
universal polynomial depending only on the singularity, if we
substitute the characteristic classes of the map. In modern language
the singularity $\eta$ is a subset of the vector space of all germs
$(\C^n,0) \to (\C^m,0)$, and its Thom polynomial is its Poincar\'e
dual in equivariant cohomology (see the group action below).
Different methods to compute certain concrete Thom polynomials
(resolution, interpolation, localization, Gr\"obner basis, etc) have
been found recently, and the enumerative and combinatorial
application of the computed Thom polynomials remains a hot area in
geometry and algebraic combinatorics.

In this paper we study the interior structure of Thom polynomials of
singularities. We carry out a program, similar in spirit to the
`Thom-Sebastiani' program of local singularity theory, by studying
Thom polynomials of product singularities. Our main example uncovers
a so far hidden strong constraint on Thom polynomials.

Using this new property we reduce the knowledge of Thom polynomials
of `contact' singularities with a given local algebra to the
knowledge of a formal power series, that we named Thom series. This
is the content of Theorem~4.1 below. Namely, under the conditions
listed there, a finite dimensional commutative local $\C$-algebra
$Q$ gives rise to a formal power series $\ts_Q$ in the variables
$d_0,d_{\pm 1},d_{\pm 2},\ldots$, such that the Thom polynomial of
any singularity $\eta:(\C^n,0)\to (\C^{n+k},0)$ with local algebra
$Q$ can be obtained from $\ts_Q$ be substituting $d_i=c_{i+k+1}$.
For example, for $Q=\C[x]/(x^3)$ we have
$\tp_Q=d_0^2+d_{-1}d_1+2d_{-2}d_2+4d_{-3}d_3+\ldots$, hence the Thom
polynomial of a singularity $\eta:(\C^n,0)\to (\C^n,0)$ with local
algebra $Q$ is $c_1^2+c_2$; the Thom polynomial of a singularity
$\eta:(\C^n,0)\to (\C^{n+1},0)$ with local algebra $Q$ is
$c_2^2+c_1c_3+2c_4$, etc.

The authors thank A. Buch, M. Kazarian, B. K\H om\H uves, and A.
N\'emethi for valuable discussions.

\section{Thom polynomials and relations between them} \label{relations}

In this paper we use cohomology with rational coefficients.

Let the Lie group $G$ act on the complex vector space $V$ with a  $G$-invariant irreducible complex subvariety $\eta\subset V$ of real codimension $l$. The Poincar\'e dual of  $\eta$ in {\sl$G$-equivariant cohomology} is called the Thom polynomial of $\eta$ and is denoted by $\tp(\eta)\in H_G^*(V)=H^*(BG)$. More generally, the kernel of the restriction ring homomorphism $H_G^*(V)\to H_G^*(V\setminus \eta)$ is called the avoiding ideal of
$\eta$ and will be denoted by $\A_\eta$. The avoiding ideal has the following remarkable property
   $$\A_\eta \cap H^i(BG)=\begin{cases}0 & \hbox{   if } i<l \\ \Q,\ \hbox{spanned by} \tp(\eta) & \hbox{   if } i=l.\end{cases}$$
For more general definitions, discussions, and the proof of this property see \cite{cr}. (The avoiding ideal is not necessarily principal; it just shares the property of homogeneous principal ideals that the lowest degree piece is 1-dimensional.)

\smallskip

Consider now the Lie groups $G$ and $G\ss=G\times H$ acting on the
vector spaces $V$ and $V\ss$ respectively. Let $j:V\to V\ss$ be a
continuous map for which $(g,h)\cdot j(v)=j(g\cdot v)$ for $g\in G,
h\in H, v\in V$. This is equivalent to requiring that $j$ is a
$G\ss$-equivariant map, where the $H$-action on $V$ is defined to be
trivial. Let $\eta\subset V$ and $\eta\ss\subset V\ss$ be invariant
subvarieties with $j^{-1}(\eta\ss)=\eta$, as in the diagram
\[\begin{matrix}
G & < & G\ss & =G\times H \\
\curvearrowright & &   \curvearrowright   &               \\
V & \mapright{j} & V\ss & \\
\cup & & \cup & \\
\eta & \mapright{} & \eta\ss. \end{matrix}\]

\begin{theorem}\label{setting} Let $p=\sum x_{i}
\otimes y_i\in \A_{\eta\ss}\vartriangleleft
H^*(BG\ss)=H^*(BG)\otimes H^*(BH)$, where $y_i$ is an additive basis
of $H^*(BH)$. Then $x_i\in \A_{\eta}\vartriangleleft H^*(BG)$ for
every $i$.
\end{theorem}

\begin{proof} The rings $H_{G\ss}^*(V)$ and $H_{G\ss}^*(V\ss)$ are
isomorphic to $H^*(BG\times BH)$ (since $V$ and $V\ss$ are contractible) and $j^*$ is
an isomorphism between them. The theorem follows from the
commutative diagram
\[\begin{matrix}
 &  H_{G\ss}^*(V) & \mapleft{j^*} & H_{G\ss}^*(V\ss) & \ni p =\sum x_i \otimes y_i\\
 & \mapdown{y_i\mapsto y_i} & & \mapdown{} & \\
H_G^*(V\setminus \eta ) \langle y_i\rangle = &
H_{G\ss}^*(V\setminus\eta) & \mapleft{} &
H_{G\ss}^*(V\ss\setminus\eta\ss),
\end{matrix}\]
where vertical arrows are restriction homomorphisms. The equality
$H_G^*(V\setminus \eta ) \langle y_i\rangle =
H_{G\ss}^*(V\setminus\eta)$ holds, since the action of $H$ on $V$ is
trivial.
\end{proof}

\bigskip

Now we specialize to examples in singularity theory. We will use
some standard notions of singularity theory, see~e.g.~\cite{agvl}.
Namely, $\E(n,m)$ will denote the infinite dimensional vector space
of holomorphic germs $(\C^n,0)\to (\C^m,0)$. On this space the
so-called right-left and the contact group $\K(n,m)$ act. Both
groups are infinite dimensional, but---in a generalized sense---they
are both homotopy equivalent to $GL(n)\times GL(m)$ \cite{rl}. (Our
base field is $\C$ in the whole paper, i.e. $GL(n)$ means
$GL(n,\C)$.)

Certain $G=\K(n,m)$ invariant subsets $\eta\subset \E(n,m)$ define
Thom polynomials
\[\tp(\eta) \in H^*(B\K(n,m))=H^*(BGL(n)\times
BGL(m))=\Q[a_1,a_2,\ldots,a_n, b_1,b_2,\ldots,b_m],\] where $a_i$
and $b_i$ are the universal Chern classes of $GL(n)$ and $GL(m)$
respectively ($\deg a_i=\deg b_i=2i$).

It is known that the Thom polynomial $\tp(\eta)$  of the {\sl
singularity} $\eta\subset \E(n,m)$  only depends on the quotient
Chern classes, i.e. on the classes $c_i$ defined by
\[1+c_1t+c_2t^2+\ldots=\frac{1+b_1t+b_2t^2+\ldots}{1+a_1t+a_2t^2+\ldots},\]
where $t$ is a formal variable. This theorem was a folklore
statement in the 60's, rediscovered by many recently. The first
written proof is probably in \cite{damonphd}---we will call this
statement the Thom-Damon theorem. Observe that the statement holds
only for contact classes and not for right-left classes in general.

\bigskip

Returning to the setting of Theorem \ref{setting}, let $V=\E(n,m)$,
$V\ss=\E(n,m+1)$ and let $j:V\to V\ss$ be defined by
\[j(f)(x_1,\ldots,x_n)=(f(x_1,\ldots,x_n),0).\]
Let $G=\K(n,m)$ be the contact group and $H=U(1)$. In the rings
\[\begin{aligned}
H^*(B\K(n,m))&= \Q[a_1,\ldots,a_n,b_1,\ldots,b_m],
\\
H^*(B(\K(n,m)\times U(1)))&=
\Q[a_1,\ldots,a_n, b_1,\ldots,b_m,y] \qquad\qquad (\deg y=2)
\\
 & =H^*(B\K(n,m))[y]
\end{aligned}
\]
we define elements $c_i$ and $c_i\ss$ by
\[\begin{aligned}
1+c_1t+c_2t^2+\ldots=&\frac{ \sum_{i=0}^m b_it^i}{\sum_{i=0}^n
a_it^i},\\
1+c\ss_1 t+c\ss_2 t^2+\ldots=&\frac{(\sum_{i=0}^m
b_it^i)(1+yt)}{\sum_{i=0}^n a_it^i}.\end{aligned}\] Here and in the
sequel we use the convention that $a_0=b_0=1$, $a_{<0}=b_{<0}=0$.

Let us choose one of the following two cases for the pair $\eta\subset
V$, $\eta\ss\subset V\ss$.
\begin{description}
\item[A]\label{one} For a finite dimensional local $\C$-algebra $Q$ let
$\eta_0\subset V$ and $\eta_0\ss\subset V\ss$ be the subsets
(`singularities') with associated local algebra $Q$. Suppose that
the closures $\eta:=\overline{\eta_0}$ and
$\eta\ss:=\overline{\eta_0\ss}$ define Thom polynomials.
\item[B]\label{two}  Let $\eta_0=\Sigma^{I}(n,m)\subset V$ and
$\eta\ss_0=\Sigma^{I}(n,m+1)\subset V\ss$ be the collection of germs
with Thom-Boardman symbol $I$ (see \cite{boardman}), and let $\eta$,
$\eta\ss$ be their closures.
\end{description}
In both cases denote the complex codimensions of $\eta$ and
$\eta\ss$ by $l$ and $l\ss$ respectively. In both cases $\eta$ and
$\eta\ss$ define Thom polynomials, and by the Thom-Damon theorem
they only depend on the `quotient' Chern classes $c_i$ and $c_i\ss$
respectively.  Indeed, in case {\bf A} this follows from the fact
that two germs are contact equivalent if and only if their local algebras are
isomorphic. In case {\bf B} this follows from the fact that
Thom-Boardman classes are contact invariant. We also have that
$j^{-1}(\eta\ss)=\eta$ so the conditions of Theorem~\ref{setting}
are satisfied. To apply Theorem~\ref{setting} to this situation, we
need the following definition and proposition.

\begin{definition} Consider infinitely many variables
$x_1,x_2,\ldots$ and use the convention of $x_0=1$, $x_{<0}=0$. Let
$p=x_{u_1}x_{u_2}\cdots x_{u_r}$ be a monomial (with the $u_j$'s not
necessarily different) in the variables $x_1,x_2,\ldots$ and let $i$
be a nonnegative integer. Let $\binom{[r]}{i}$ denote the set of
$i$-element subsets of $\{1,\ldots,r\}$, and for such a subset $I$,
$\delta_I$ will denote its characteristic function. We define the
lowering operator $\flat[i]$ by
\[p\fl{i}(x_1,x_2,\ldots)=\sum_{I\in\binom{[r]}{i}} x_{u_1-\delta_I(1)}
x_{u_2-\delta_I(2)}\cdots x_{u_r-\delta_I(r)},\] and we extend the
definition to all polynomials linearly.
\end{definition}

\noindent For example $(x_1x_2x_5+x_8+x_4^2)\fl{2}=x_1x_5+x_2x_4+x_1^2x_4+x_3^2$.

\begin{proposition}\label{flat} Let $p(x_1,x_2,\ldots)$ be a degree $l \ss$ polynomial for which
\[p(c\ss_1, c\ss_2,\ldots)=\sum_i p_i(c_1,c_2,\ldots)y^{l\ss-i}.\] Then $p_i=p\fl{i}$.
\end{proposition}

\begin{proof} Since $1+c\ss_1 t+c\ss_2 t^2+\ldots=(1+c_1 t+c_2t^2+\ldots)(1+yt)$, the proposition follows from the definition of $\flat[i]$. \end{proof}

Now we are able to state the connection between the Thom polynomials $\tp({\eta\ss})$ and $\tp(\eta)$.

\begin{theorem}\label{tpflat} Let $\tp(\eta\ss)=p(c\ss_1,c\ss_2,\ldots)$ and $\tp(\eta)=q(c_1,c_2,\ldots)$. Then
\[p\fl{i}=\begin{cases} 0 & \text{if } i>l\ss-l \\ k_\eta \cdot q, & \text{if } i=l\ss-l  \end{cases}\]
where $k_\eta\in\Q$.
\end{theorem}

\begin{proof} The polynomial
$p\fl{i}$ must belong to the avoiding ideal $\A_\eta$ according to
Theorem \ref{setting} and Proposition~\ref{flat}. According to
\cite[Sect. 2.3]{cr}---recalled at the beginning of this
section---the lowest degree elements of the avoiding ideal
$\A(\eta)$ have degree $l$, and are constant multiples of the Thom
polynomial $\tp(\eta)$.
\end{proof}

In fact $k_\eta\in\Z$ since the cohomology of $BGL(n)$ has no
torsion so we can work here with integer coefficients.

\begin{corollary}\label{numbfactors} The number of factors in any term of $\tp({\eta\ss})(c\ss_1,c\ss_2,\ldots)$ is at most $l\ss-l$. \end{corollary}

\begin{proof}
Let the longest terms have $i>l\ss-l$ factors, and let $T$ be the
lexicographically largest term with $i$ factors. Then  $T\fl{i}$ has
one term, and that term will not be among the terms of $S\fl{i}$,
for the other terms $S$. Hence $\tp({\eta\ss})\fl{i}\not=0$, which
contradicts to Theorem \ref{tpflat}.
\end{proof}

We can assume that all terms of $\tp(\eta\ss)$ have exactly $l\ss-l$
factors, since we can multiply the shorter terms by an appropriate
$c\ss_0$-power. The benefit is that the $\flat[l\ss-l]$ operation is
particularly simple for polynomials with terms of exactly $l\ss-l$
factors: it decreases all indices by 1. Hence it is worth shifting
the indices in the notation: let $d_i=c_{i+(m-n+1)}$ and
$d\ss_i=c\ss_{i+(m-n+2)}$. Then, in the cases when $k_\eta=1$, we
have that the Thom polynomial of $\eta$ (written in the
$d$-variables) is obtained from the Thom polynomial of $\eta\ss$
(written in the $d\ss$-variables) by deleting the $\sharp$'s. For
example we have
\[\begin{aligned}
\hbox{for } m=n+1 & \phantom{xxx} & \tp({A_2})=& c\ss_2 c\ss_2+c\ss_1c\ss_3+2c\ss_4&=d\ss_0d\ss_0+d\ss_{-1}d\ss_1+2d\ss_{-2}d\ss_2, \\
\hbox{for  } m=n & \phantom{xxx} & \tp({A_2})=& c_1c_1+\ \ \
c_2&=d_0d_0+d_{-1}d_1+2d_{-2}d_2
\end{aligned}\]
cf. Section \ref{thomseries}. Observe that $d_{-2}d_2$ is 0 in the
second line, but not 0 in the first line. Hence from
$\tp({\eta\ss})$ we can compute $\tp({\eta})$ but not vice versa in
general.

\section{The constant $k_\eta$.}

In this section we study the constant $k_\eta$ occurring in Theorem
\ref{tpflat}. We give conditions for the non-vanishing of it, based
on the observation that $\tp(\eta)$ restricted to the smooth part
$\eta_0$ of $\eta$ is the equivariant normal Euler class of
$\eta_0$. This property was crucial in all the Thom polynomial
calculations using the method of restriction equations in
\cite{rrtp}.

Let $V=\E(n,m)$, $V\ss=\E(n,m+1)$, $G=\K(n,m)$, $H=U(1)$, $j:V\to
V\ss$, $\eta_0\subset\eta$, $\eta\ss_0\subset\eta\ss$ (case {\bf A}
or {\bf B}), $\tp(\eta)=q(c_1,c_2,\ldots)$, and
$\tp(\eta\ss)=p(c\ss_1,c\ss_2,\ldots)$ be defined as as above; and
for $X\subset Y$ let $\nu_X$ denote the normal bundle of $X$.
Observe that through $j$ we can identify the bundle $\nu_{\eta_0}$
as a subbundle of $\nu_{\eta\ss_0}|_{j(\eta_0)}$, hence for the
$G\times U(1)$-equivariant Euler classes we have
\[e\Big(\nu_{\eta\ss_0}\Big|_{j(\eta_0)}\Big)=e(\nu_{\eta_0}) \cdot  e(M),\]
where $M$ is the complementary subbundle of rank $l\ss-l$ (Replacing $G$ with
the subgroup
$U(n)\times U(m)$---it doesn't change the equivariant cohomology---we can
assume that $M$ is $G$-equivariant). Let
$r:H_G^*(V)\to H_G^*(\eta_0)$, and $r\ss:H_{G\times U(1)}^*(V\ss)\to
H_{G\times U(1)}^*(\eta\ss_0)$ be the restriction homomorphisms.
Then---since the cohomology class of a submanifold restricted to the
submanifold itself is the Euler class of the normal bundle of the
submanifold---we obtain that
\begin{equation}\label{a1}p(j^*r\ss(c\ss_1),j^*r\ss(c\ss_2),\ldots)=q(r(c_1),r(c_2),\ldots)\cdot e(M).\end{equation}
Using the identity
  \[1+c\ss_1t+c\ss_2t^2+\ldots=(1+yt)(1+c_1t+c_2t^2+\ldots)\]
and Corollary \ref{numbfactors}, the left hand side can be written
as
  \[y^{l\ss-l}\cdot p\fl{l\ss-l}\big(r(c_1),r(c_2),\ldots\big)+\text{terms of lower }y\text{-power}.\]
According to Theorem \ref{tpflat} it is further equal to
  \[y^{l\ss-l}\cdot k_\eta \cdot q(r(c_1),r(c_2),\ldots)+\text{terms of lower }y\text{-power}.\]
Thus from equation (\ref{a1}) we obtain that if
the normal Euler class $e(\nu_{\eta_0})=q(r(c_1),r(c_2),\ldots)$ is
not 0, then $k_\eta$ is the coefficient of $y^{l\ss-l}$ in $e(M)$.
Restricting the action of $G\times U(1)$ to $U(1)$ we obtain the
following equivalent statement.

\begin{theorem}\label{cons} Suppose $e(\nu_{\eta_0})\not=0$. Then the $U(1)$-equivariant Euler class
$e_{U(1)}(M)=k_\eta \cdot y^{l\ss-l} \in H^*(BU(1))=\Q[y]$.
\end{theorem}

The points of $\eta_0$ are $U(1)$-invariant, hence the coefficient
$k_\eta$ is the product of the weights of the $U(1)$-action on a
fiber of $M$. When $\eta_0$ is the contact singularity
corresponding to the local algebra $Q$ the fiber of $M$ is identified with
the maximal ideal of the local algebra $Q$ with scalar $U(1)$-action
(see \cite{rrtp}), that is, in this case $e_{U(1)}(M)=1$.

It also follows from results is \cite{rrtp} that the vanishing of
$e(\nu_{\eta_0})$ only depends on the local algebra of $\eta_0$, not
on the particular dimensions $n$ and $m$. Also, the Euler class
$e(\nu_{\eta_0})$ is not zero for any of the {\sl simple}
singularities considered in that paper---e.g. singularity types
$A_i$, $I_{a,b}$, $III_{a,b}$, and more. (In fact for all
singularities with known Thom polynomial the Euler class
$e(\nu_{\eta_0})$ is not zero.) Hence for all these singularities we
proved that Theorem~\ref{tpflat} holds with constant $k_\eta=1$.

\section{Thom series} \label{thomseries}

Now we use the $d_i=c_{i+(m-n+1)}$ substitution to write the Thom polynomials of natural infinite sequences of singularities in a concise form.

\subsection{Contact singularities.}\label{contactsing}

Let $Q$ be a local algebra of a singularity. We will need three
integer invariants of $Q$ as follows: (i) $\delta=\delta(Q)$ is the
complex dimension of $Q$, (ii) the {\em defect} $d=d(Q)$ of $Q$ is
defined to be the minimal value of $b-a$ if $Q$ can be presented
with $a$ generators and $b$ relations; (iii) the definition of the
third invariant $\gamma(Q)$ is more subtle, see \cite[\S6]{mather6}.
The existence of a {\em stable} singularity $(\C^n,0)\to (\C^m,0)$
with local algebra $Q$ is equivalent to the conditions $m-n\geq d$,
$(m-n)(\delta-1)+\gamma\leq n$. Under these conditions the
codimension of the $\K(n,m)$ orbit of a germ
$\eta:(\C^n,0)\to(\C^m,0)$ with local algebra $Q$ is
$(m-n)(\delta-1)+\gamma$.

Now we can apply Theorem~\ref{tpflat} and
Corollary~\ref{numbfactors} to the series of contact singularities
with local algebra $Q$. Observe that Theorem \ref{cons} and the
discussion after computes the value $k_\eta=1$. Hence we obtain the
following theorem.
\begin{theorem} \label{q} Let $Q$ be a local algebra of singularities
defining Thom polynomials. Assume that the normal Euler classes of
these singularities are not 0. Then associated with $Q$ there is a
formal power series ({\em Thom series}) $\ts_Q$ in the variables
$\{d_i|i\in \Z\}$, of degree $\gamma(Q)-\delta(Q)+1$, such that all
of its terms have $\delta(Q)-1$ factors, and the Thom polynomial of
$\eta:(\C^n,0)\to (\C^m,0)$ with local algebra $Q$ is obtained by
the substitution $d_i=c_{i+(m-n+1)}$. \qed
\end{theorem}

Even though there are powerful methods by now to compute individual Thom polynomials (i.e.~finite initial sums of the $\ts$), finding closed formulas for these Thom series remains a subtle problem. Here are some examples.

\begin{description}
\item[\bf A$_0$] $Q=\C$ ({\sl embedding}). Here $\delta=1$, $\gamma=0$, and
\[\ts=1.\]
\item[\bf A$_1$] $Q=\C[x]/(x^2)$ (e.g. {\sl fold, Whitney umbrella}). Here $\delta=2$, $\gamma=$1, and
\[\ts=d_0.\]
\item[\bf A$_2$] $Q=\C[x]/(x^3)$ (e.g. {\sl cusp}). Here $\delta=3$, $\gamma=2$, and \cite{rongaij}
\[\ts=d_0^2+d_{-1}d_1+2d_{-2}d_2+4d_{-3}d_3+8d_{-4}d_4+\ldots \]
\item[\bf A$_3$] $Q=\C[x]/(x^4)$. Here $\delta=4$, $\gamma=3$, and \cite[Thm.4.2]{a3}
\[\ts=\sum_{i=0}^\infty 2^i
d_{-i}d_0d_i+\frac{1}{3}\sum_{i=0}^\infty \sum_{j=0}^\infty 2^i3^j
d_{-i}d_{-j}d_{i+j}+\frac{1}{2} \sum_{i=0}^\infty \sum_{j=0}^\infty
a_{i,j}d_{-i-j}d_id_j,\] where $a_{i,j}$ is defined by the formal
power series
\[\sum_{i,j}
a_{i,j}u^iv^j=\frac{u\frac{1-u}{1-3u}+v\frac{1-v}{1-3v}}{1-u-v}.\]
\item[\bf I$_{2,2}$] $Q=\C[x,y]/(xy,x^2+y^2)$. Here $\delta=4$, $\gamma=4$, and
\[\ts=\sum_{i=1}^{\infty} 2^{i-2}d_{-i}d_1d_i-\sum_{i=1}^\infty
2^{i-1}d_{-i}d_0d_{i+1}+\frac{1}{2}\sum_{i=1}^\infty\sum_{j=1}^\infty
\binom{i+j-2}{i-1} d_{-i-j+1}d_id_j.\]
\end{description}

Strictly speaking we have not proved the Thom series of $I_{2,2}$,
just obtained overwhelming computer evidence for it. For an accurate
proof (using the method of restriction equations from \cite{cr}) we
would need to manipulate non-trivial resultant identities similar to
those in \cite{a3}.

Recently P. Pragacz in \cite{pragacz} used Schur functions and other
powerful tools of symmetric functions in order to calculate the Thom
polynomials of certain contact class singularities $\eta:\C^n\to
\C^m$ for any $n$ and $m$. In particular he calculates the Thom
polynomials of the singularities of the list above in terms of Schur
functions. We do not see how to prove directly that the two results
agree. We were informed by A. Szenes and G. B\'erczi that they
calculated Thom series for $A_i$ for some values of $i>3$.

\begin{remark}
A remarkable property of all the computed Thom series is that they
have positive coefficients when written in the basis of Schur
polynomials instead of the basis of Chern monomials. This was
recently proved in \cite{pragweber}.
\end{remark}

\subsection{Thom-Boardman classes.}

The Thom-Boardman singularity $\Sigma^I$ is the closure of a contact
singularity if the relative dimension $m-n$ is large. In this case
the constant $k_\eta$ is 1, due to the results of
Section~\ref{contactsing}. We conjecture that the constant $k_\eta$
is 1 for Thom-Boardman singularities with arbitrary relative
dimension $m-n$, too. However it is only known for the cases
$\Sigma^i$ and $\Sigma^{i,j}$.

\begin{theorem} \label{sigmai} For every $r$ there is a formal power series
$\ts_{\Sigma^r}$ in the variables $\{d_i|i\in \Z\}$, of degree
$r(r-1)$, such that all of its terms have $r$ factors, and the Thom
polynomial of $\Sigma^r(n,m)$ is obtained by the substitution
$d_i=c_{i+(m-n+1)}$.
\end{theorem}

This is the statement we obtain by applying Theorem~\ref{tpflat} to
the corank $r$ singularities (case {\bf B} above with $I=(r)$). It
is not new though, since all these polynomials are known explicitly
(Giambelli-Thom-Porteous):
\[\tp(\Sigma^r)=\det(d_{r-1+j-i})_{i,j=1,\ldots, r},\]
that is, the Thom series in this case is finite. Also observe that
to prove Theorem~\ref{sigmai} we did not need the infinite
dimensional spaces $\E(n,m)$, we could have started with
$\Hom(\C^n,\C^m)$.

The codimension of the set of germs $\C^n\to \C^m$ with
Thom-Boardman symbol $(i,j)$ is $(i+m-n)i+j((i+m-n)(2i-j+1)-2(i-j))/2$.
There are algorithms to calculate  the Thom polynomials of second
order Thom-Boardman singularities see \cite{rongaij}, \cite{kazaij},
\cite{kf}. From \cite{kf} it also follows that $k_\eta=1$ for
$\Sigma^{i,j}$ singularities. Hence we have the following result
about the structure of Thom polynomials of singularities of
Thom-Boardman symbol $(i,j)$.

\begin{theorem} \label{sigmaij}
For every $i\geq j$ there is a formal power series
$\ts_{\Sigma^{i,j}}$ in the variables $\{d_i|i\in \Z\}$ of degree
$i(i-1)+j(2i^2-ij-3i+3j-1)/2$, such that all of its terms have
$i(j+1)-\binom{j}{2}$ factors, and the Thom polynomial of
$\Sigma^{i,j}(n,m)$ is obtained by the substitution
$d_i=c_{i+(m-n+1)}$.
\end{theorem}

In general closed formulas for the  $\Sigma^{i,j}$ Thom series are not known, with the
following exceptions. The Thom series of $\Sigma^{i,1}$ is
explicitly computed in Theorem 4.8 of \cite{kf}. A closed formula for the `initial term'
of the Thom series of $\Sigma^{i,j}$ is calculated in 4.6 of
\cite{kf}. Here `initial term' refers to the terms not containing a
$d_u$ factor with $u\leq -i+1$.

\section{Other stabilizations: product singularities}

Another way of looking at the results in Section \ref{relations} is
that we related the Thom polynomial of a singularity
$\eta:(\C^n,0)\to (\C^m,0)$ to the Thom polynomial of $\eta\times
\xi$, where $\xi$ is the map $(\C^0,0)\to (\C^1,0)$. We can,
however, choose other $\xi$ maps (with some care), and relate the
Thom polynomial of $\eta$ with the Thom polynomial of $\eta\times
\xi$.

Let us fix a germ $\xi:(\C^a,0)\to (\C^b,0)$ whose right-left
symmetry group is $H$ (or $H$ is a subgroup of the right-left
symmetry group), with representations $\lambda_0$, $\lambda_1$ on
the source and the target spaces respectively. Let $c(\xi)$ be the
formal quotient Chern class $c(\lambda_1)/c(\lambda_0)\in H^*(BH)$.
Consider the map $j:\E(n,m) \to \E(n+a,m+b)$, $f\mapsto f\times \xi$.
If $\eta\subset \E(n,m)$ is the closure of a $\K(n,m)$ orbit of
codimension $l$ and
\begin{itemize}
\item $\eta\ss:=\overline{\K(n+a,p+b)\cdot j(\eta)}$ defines a Thom polynomial of degree $l\ss$, and
\item $j^{-1}(\eta\ss)=\eta$,
\end{itemize}
then Theorem~\ref{setting} gives the following relation between the
Thom polynomials of $\eta$ and $j(\eta)$.

\begin{theorem} \label{product} Let the classes $c\ss_i$ be defined by
$1+c\ss_1+c\ss_2+\ldots:=$ $(1+c_1+$ $c_2+\ldots)\cdot c(\xi)$.
Suppose that substituting these classes into $\tp(\eta\ss)$ results
in $\sum x_i \cdot y_i$. That is, let
\[\tp(\eta\ss)\Big( (1+c_1+c_2+\ldots)\cdot c(\xi)\Big) = \sum x_i \cdot y_i\qquad \in
\Q[c_1,c_2,\ldots]\langle y_i \rangle,\] where
$x_i\in\Q[c_1,c_2,\ldots]$, and $y_i$ is an additive basis of
$H^*(BH)$. Then $x_i=0$ for $\deg x_i<l$, and $x_i=constant\cdot
\tp(\eta)$ if $\deg x_i=l$.
\end{theorem}

The case studied in Section \ref{relations} is recovered as
$\xi:(\C^0,0)\to (\C^1,0)$, $H=U(1)$, $\lambda_1=$the standard
representation, $\lambda_0=$the 0-dimensional representation,
$c(\xi)=1+y\in H^*(BH)$.

The stability of Thom polynomials of contact singularities under
trivial unfolding corresponds to the following case:
$\xi:(\C^1,0)\to (\C^1,0)$, $x\mapsto x$, $H=U(1)$,
$\lambda_0=\lambda_1=$the standard representation.

Unfortunately, it is hard to find other $\xi$'s, for which one can
easily check the two conditions above. If $\xi$ is a complicated
germ, then the set $\eta\ss$ is usually within the realm of moduli
of singularities, where the question of {\sl what} defines a Thom
polynomial is very hard (it reduces to the computation of Kazarian's
spectral sequence \cite{kaza}, \cite[Sect.10]{cr}). Hence we show an
example with a simple $\xi$.

\begin{example}\label{masik} \rm
Let $\xi:(\C,0)\to (\C,0)$, $x\mapsto x^2$; and let $H=U(1)$ act by
$\lambda_0=\rho$, $\lambda_1=\rho^{\otimes 2}$, where $\rho$ is the
standard representation of $U(1)$. Let $\eta\subset \E(n,n)$ be the
closure of the 2-codimensional $\K$-orbit corresponding to the local
algebra $\C[x]/(x^3)$, i.e. the orbit of
$(x_1,x_2,\ldots,x_n)\mapsto$ $(x_1^3,x_2,$ $\ldots,x_n)$. Then
$\eta\ss=j(\eta)$ is the closure of the 7-codimensional
$\K(n+1,n+1)$ orbit corresponding to the local algebra
$\C[x,y]/(x^3,y^2)$. The two conditions above hold and we obtain,
that if
\begin{equation}\tp(\eta\ss)\Big(
(1+c_1+c_2+\ldots)(1+y-y^2+y^3-\ldots)\Big)=x_0y^7+x_1y^6+\ldots+x_7y^0,
\label{itt}\end{equation} then we must have $x_0=x_1=0$ and
$x_2=constant\cdot \tp(\eta)$. Indeed, these polynomials are known
explicitly,
$\tp(\eta\ss)=2(c\ss_1{c\ss_2}^3-{c\ss_1}^2c\ss_2c\ss_3+{c\ss_2}^2c\ss_3+c\ss_1{c\ss_3}^2-2c\ss_1c\ss_2c\ss_4+c\ss_3c\ss_4-c\ss_2c\ss_5)$
\cite{por}, and thus the left hand side of (\ref{itt}) becomes
$2\Big((c_1+y)(c_2+c_1y-y^2)^3-\ldots\Big)=
4(c_1^2+c_2)y^5+2(c_1c_2-c_1^3+10c_3)y^4+\ldots$. This is consistent
with the fact that $\tp(\eta)=c_1^2+c_2$ (and the other $y$
coefficients also belong to $\A_\eta$).
\end{example}

\smallskip
Notice that Theorem \ref{product} says that under the specified conditions we can calculate  $\tp(\eta)$ and $\tp(\xi)$ if we know $\tp(\eta\times \xi)$ but not in the other direction. Nevertheless it gives strong restrictions on the form of $\tp(\eta\times \xi)$.

The analogous `Thom-Sebastiani' approach to finding relations
between different Thom polynomials is more promising in finite
dimensional settings, where we know the existence of Thom
polynomials, and their properties are relevant in algebraic
combinatorics. We plan to study the Thom polynomials for quiver
representations from this perspective in the future.

\bigskip

\begin{remark}
There are other natural infinite series of pairs of Thom polynomials
where Theorem~\ref{tpflat} applies. For instance, we can consider
the Thom polynomials of the orbits of $S^2(\C^n)$ and
$S^2(\C^{n+1})$ with actions of $GL(n)$ and $GL(n)\times U(1)$.
Their Thom polynomial theory is worked out, see \cite{jlt}, also
\cite[Sect. 5]{cr}. The constant $k_\eta$ in this case, however,
turns out to be different from 1.
\end{remark}

\bibliography{dstab}

\bibliographystyle{alpha}

\end{document}